\documentclass[11pt]{article}
\usepackage{amsmath}
\usepackage{amssymb}
\usepackage{epsf}
\usepackage{amscd}
\usepackage[all]{xy}
\usepackage[francais,english]{babel}

\newtheorem{theorem}{Theorem}[section]
\newtheorem{propo}[theorem]{Proposition}
\newtheorem{coro}[theorem]{Corollary}
\newtheorem{lemme}[theorem]{Lemma}
\newtheorem{defn}[theorem]{Definition}
\newtheorem{propdefe}[theorem]{Definition-Proposition}
\newtheorem{rquee}[theorem]{Remark}
\newtheorem{csqcee}[theorem]{Consequence}
\newtheorem{exple}[theorem]{Example}
\newtheorem{proptee}{Property}
\newtheorem{notationse}[theorem]{Notation}

\newenvironment{theo}{\begin{theorem} \sf}{\end{theorem}}
\newenvironment{prop}{\begin{propo} \sf}{\end{propo}}
\newenvironment{cor}{\begin{coro} \sf}{\end{coro}}
\newenvironment{lemma}{\begin{lemme} \sf}{\end{lemme}}
\newenvironment{dfn}{\begin{defn} }{\end{defn}}
\newenvironment{proof}{{\flushleft{\bf Proof: }} \rm}{\Qed \\}
\newenvironment{prooflemma}{{\flushleft{\bf Proof of Lemma: }} \rm}{\qed \\}
\newenvironment{defprop}{\begin{propdefe} }{\end{propdefe}}

\newenvironment{rmk}{\begin{rquee} \normalfont}{\end{rquee}}
\newenvironment{ex}{\begin{exple} \normalfont}{\end{exple}}

\def\qed{\hfill \mbox{$\square$}}
\def\Qed{\hfill \mbox{$\blacksquare$}}
\def\ot{\otimes}
\def\hot{\widehat{\otimes}}

\def\id{\mathrm{id}}
\def\uot{\mbox{\b{$\otimes$}}}
\def\oot{\bar{\otimes}}
\def\qot{\uot _H}

\newcommand{\hb}{\mathrm{H}_b^*(H,H)}

\newcommand{\Ext}{\mathrm{Ext}}

\newcommand{\Hh}{\mathrm{H}}

\newcommand{\Hom}{\mathrm{Hom}}

\newcommand{\hopf}{\mathcal{H}}

\newcommand{\env}{\mathcal{E}}

\newcommand{\ext}{\mathcal{E}\!xt _{\hopf}}
\newcommand{\e}{{\bf E}}
\newcommand{\f}{{\bf F}}
\newcommand{\g}{{\bf G}}

\numberwithin{equation}{section}

\title{\textsc {Injective Hopf bimodules, cohomologies of infinite
  dimensional Hopf algebras and graded-commutativity of the Yoneda
  product.}}

\author{Rachel Taillefer\thanks{Supported by a Lavoisier grant from
    the French Ministry of Foreign Affairs.}\\ \\ \footnotesize{ 
Department of Mathematics and Computer
    Science,}\\ \footnotesize{University of Leicester,}\\
  \footnotesize{Leicester LE1 7RH,} \\ \footnotesize{United
    Kingdom.}\\ \footnotesize{E-mail: R.Taillefer@mcs.le.ac.uk.}}

\begin{document}

\maketitle


\selectlanguage{english}

\begin{abstract} 

We prove that the category of Hopf bimodules over any Hopf algebra  has enough injectives, which enables us to extend
some results on the unification of Hopf bimodule cohomologies of \cite{Ta,Ta2} to the infinite dimensional case. We
also prove that the cup-product defined on these cohomologies is
graded-commutative.

\end{abstract}


\section{Introduction}

The object of this paper is to extend the results in \cite{Ta,Ta2},  on cohomologies for
Hopf algebras, to
a more general and more useful context, and
to prove some properties of the objects introduced there. 

For a finite dimensional Hopf algebra $H$, we
identified various cohomologies with the $\Ext^*$ functor over an `enveloping' associative
algebra of $H$ introduced by C. Cibils and M. Rosso
(\cite{CR}). In this way, the category of Hopf bimodules over $H$ and
the category of modules over this enveloping algebra can be
identified, and we are then able to use all the usual properties of
the $\Ext^*$ functor for a module category. When $H$ is infinite
dimensional however, such an `enveloping' algebra
does not exist. Nevertheless, many interesting examples of Hopf algebras are in fact infinite
dimensional (\emph{e.g.} enveloping algebras of Lie algebras), and it
is important that such cases be considered. 
In fact, the identification of
the cohomologies with an $\Ext^*$ functor (still for
the category of Hopf bimodules over $H$, but this category is not equivalent to
a module category in this case) does still hold: this
generalisation is the object of Sections \ref{cohomologies} and
\ref{enoughinjectives}.  The obstacle to this generalisation was the
existence of enough projectives in the category of Hopf bimodules over
$H$: it
is not known whether any Hopf bimodule is a quotient of a projective
Hopf bimodule. However, this problem can be circumvented: the proofs
of the identification (which use universal properties of the $\Ext^*$
functor) can be adapted so that the obstacle becomes the
existence of enough \emph{injectives} in the category of Hopf
bimodules (\emph{see} Theorem \ref{unificationtheorem}), and this
obstacle can be removed: we prove in
Section \ref{enoughinjectives} that there are indeed
enough injectives in the category of Hopf bimodules, that is that any
Hopf bimodule can be embedded in an injective Hopf bimodule. As a
result, the identification of the cohomologies will hold for any Hopf
algebra, that is each of the cohomologies we are considering is
isomorphic to the $\Ext^*$ functor for the category of Hopf bimodules.

We also defined in \cite{Ta,Ta2} a cup-product on some of these
cohomologies, and proved that it corresponds to the Yoneda product of
extensions \emph{via} this identification. (This correspondence was
proved when the Hopf algebra is finite
dimensional, since we needed the identification, but the proofs are valid for infinite dimensional Hopf
algebras). We are in fact interested
in the algebraic structure of one of these cohomologies, $\hb,$ which is
without coefficients: we would like to know whether this cohomology is
a Gerstenhaber algebra. Gerstenhaber algebras are graded algebras,
which are graded-commutative, and which are endowed with a graded Lie
product (for a different grading) compatible with the
cup-product. This graded Lie bracket appears in many cohomologies
adapted to the study of deformations of various structures; it
describes the obstructions to deforming
 the structures for which these cohomologies are
defined. It also enabled M. Gerstenhaber and S.D. Schack in \cite{GS92}
to give a
short proof of the Hochschild-Kostant-Rosenberg theorem.

In
this paper, we prove that the cup-product is graded-commutative, using
techniques of
S. Schwede (\emph{cf.} \cite{S}), which tends to
encourage the possibility of the cohomology being a Gerstenhaber
algebra. We also give a candidate for a graded Lie bracket, which is graded-anticommutative.

This paper is organised as follows: in the first section we recall the
definitions of the cohomologies and give some of the proofs for the unification (using
injective Hopf bimodules) in the finite dimensional case. In the next
section we  prove
that every Hopf bimodule over any Hopf algebra is a sub-Hopf bimodule of an injective Hopf bimodule, and that as a
consequence the unification is true for an infinite dimensional Hopf
algebra. We then consider
the algebraic structure of $\hb$,
first proving that the cup-product defined in \cite{Ta,Ta2} is
graded-commutative in Section \ref{commut}, and
finally in the last section we give a candidate for a graded Lie
bracket.

{\bf Acknowledgement:}  I would like to thank Marco A. Farinati for his
comments on the proof of Theorem \ref{injectives}.

\section{Unification of cohomologies associated to a finite
dimensional Hopf algebra}\label{cohomologies}

Let us define first of all the cohomologies we are
interested in, and briefly give the proofs of the identification. In
all the following, $k$ is a commutative field.

\subsection{Preliminaries}

Let $(H,\mu,\eta,\Delta,\varepsilon,S)$ be a Hopf algebra.  We shall define some Hopf bimodule structures on the
  tensor product of Hopf bimodules; $\ot$ denotes the tensor product
  over $k.$

\begin{defprop} Let $M$ and $N$ be Hopf bimodules over $H,$ and consider the vector
space $M\ot N.$ This can be endowed with two different Hopf bimodule
structures. 

The first one we shall denote by $M \uot N;$ its actions are the
regular ones and its coactions are codiagonal: \begin{equation*}\begin{array}{rrcl}
\mu _L : & H \otimes M \uot N & \longrightarrow & M
\uot N;\\ & h \otimes m \otimes n & \mapsto & hm
\otimes n \\ \\
 \mu _R : & M \uot N \ot H & \longrightarrow & M \uot N; \\
 &  m \otimes n \ot h & \mapsto & m \otimes nh \\ \\
 \delta _L: & M \uot N  & \longrightarrow & H \otimes M
 \uot N;\\& m \otimes n & \mapsto &  m_{(-1)} n_{(-1)}
 \otimes m_{(0)} \otimes n_{(0)}\\ \\  \delta _R: & M \uot N  & \longrightarrow & M \uot N \ot H; \\ 
&  m \otimes n & \mapsto &  m_{(0)} \otimes n_{(0)} \ot m_{(1)} n_{(1)} .
\end{array} \end{equation*}

Dually, we denote by $M \oot N$ the Hopf bimodule whose actions are
diagonal and whose coactions are regular:  \begin{equation*}\begin{array}{rrcl}
\mu _L : & H \otimes M \bar{\otimes} N & \longrightarrow & M
 \bar{\otimes} N;\\ & h \otimes m \otimes n & \mapsto &  h^{(1)}
m \otimes h^{(2)} n  \\ \\
\mu _R : & M \bar{\otimes} N \ot H & \longrightarrow & M \bar{\otimes} N;\\ 
&  m \otimes n \ot h & \mapsto & m h^{(1)}  \otimes n h^{(2)} \\ \\
\delta _L : & M \bar{\otimes} N & \longrightarrow & H \otimes M
\bar{\otimes} N;\\ & m\otimes n & \mapsto &  m_{(-1)} \otimes
m_{(0)} \otimes n\\ \\ \delta _R : & M \bar{\otimes} N & \longrightarrow &  M \bar{\otimes} N \ot H;\\
& m\otimes n & \mapsto &  m \otimes n_{(0)} \otimes n_{(1)}.
\end{array}\end{equation*} \end{defprop}

\begin{rmk}\label{inj} Suppose $M$ is an $H$-bimodule. Then, with the
  structures described
  above, the space $H\oot M \oot H$ is a well-defined Hopf bimodule. \end{rmk}

We shall now recall the properties of the bar and cobar resolutions of
Hopf bimodules (\emph{see} for instance \cite{Sh-St} or \cite{Ta}):

\begin{prop} Let $M$ and $N$ be Hopf bimodules. We view the spaces in
the right module (\emph{resp.} bimodule) bar resolution of $M$ and
the spaces in the left comodule (resp. bicomodule) cobar resolution of $N$ as
Hopf bimodules, with $Bar_q (M)=M \uot H^{\uot q +1}$
(\emph{resp.} $B_q (M)= H^{\uot q +1} \uot M \uot
H^{\uot q +1}$) and $Cob^p (N)=  H^{\oot p +1} \oot N$
(\emph{resp.} $C^p (N)=  H^{\oot p +1} \oot N \oot
H^{\oot p +1}$).

Then these resolutions are Hopf bimodule resolutions. Furthermore,
$Bar_\bullet (H)$ and 
$B_\bullet(M)$ split as  sequences of bicomodules, and dually
$Cob^\bullet (H)$ and $C^\bullet (N)$ split as sequences of
bimodules. \end{prop}

These resolutions have more properties which we shall need. To
describe them, we need some definitions:

\begin{dfn} {\rm (\cite{Sh-St})} Let $\Hom_{\hopf}(M,N)$ denote the space
  of Hopf bimodule morphisms from $M$ to $N.$ A  Hopf bimodule $N$ is called a
\emph{relative injective} if the functor $\Hom_{\hopf}(-,N)$ takes exact
sequences of Hopf bimodules that split as sequences of bimodules to
exact sequences of $k-$vector spaces.
\end{dfn}

\begin{ex}  An injective Hopf bimodule is a relative injective. \end{ex}

\begin{ex}  {\rm (\cite{Sh-St})}\label{exrelproje} If $V$ is a bimodule, then the
Hopf bimodule $H {\oot} V {\oot} H $ is a relative
injective.
\end{ex}

\begin{dfn} A resolution of a Hopf bimodule is a \emph{relative
injective resolution} if all its terms are relative injectives, and if it
splits as a sequence of bimodules.
\end{dfn}

Relative injective resolutions have properties which are similar to
those of injective resolutions, and we shall use:

\begin{prop}\label{equivse}{\rm (\emph{cf.} \cite{Sh-St}
Proposition 10.5.3)}  Two relative
injective resolutions are homotopy equivalent as Hopf bimodule
complexes.\end{prop}

We can define relative projectives and relative projective
resolutions in a dual way. Then:

\begin{prop} The bar resolutions  $Bar_\bullet (H)$ and 
$B_\bullet(M)$ are relative projective resolutions, the cobar
resolutions $Cob^\bullet (H)$ and $C^\bullet (N)$ are relative
injective resolutions. \end{prop}

We now have the background necessary to define and unify the
cohomologies.

\subsection{The cohomologies}

We are interested in three cohomologies: two were defined by
M. Gerstenhaber and S.D. Schack in \cite{GS90} and \cite{GS92}, the
third by C. Ospel in his thesis \cite{Os}. They are denoted
 by $\Hh_{GS}^*$, $\Hh_b^*$ and $\Hh_{H4}^*$ respectively, and
are defined as
follows:

\begin{dfn} The three cohomologies are all defined using the bar and
  cobar resolutions. Let $M$ and $N$ be Hopf bimodules: 

\begin{enumerate} 
\item $\Hh_{GS}^*(M,N)$ is the cohomology of the double complex
  $\Hom_{\hopf}(B_\bullet (M),  C^\bullet (N));$  at point
  $(p,q)$ the entry of the double complex is $\Hom_{\hopf}(B_q (M),  C^p
  (N)),$ the space of Hopf bimodule maps from $B_q (M)$ to $ C^p
  (N),$ the vertical differential is composition with the bar
  differential, and the horizontal differential is composition with
  the cobar differential. 

\item  $\hb$ is the cohomology of the double complex
  $\Hom_{\hopf}(Bar_\bullet (H),  Cob^\bullet (H)).$ This cohomology does
  not have coefficients.

\item  $\Hh_{H4}^*$  is the cohomology of the double complex
  $\Hom_{\hopf}(Bar_\bullet (M),  Cob^\bullet (N)).$
 \end{enumerate}\end{dfn}

\begin{rmk} It is obvious from the definitions that the third
  cohomology $\Hh_{H4}^*$ generalises the second $\Hh_b^*$. In fact, $\Hh_{GS}^*$ also
  generalises $\Hh_b^*$, and is furthermore isomorphic to
  $\Hh_{H4}^*$. To prove this, we have identified each one with the
  $\Ext^*$ functor
  over an `enveloping' algebra $X$ of $H,$ defined by C. Cibils and
  M. Rosso:\end{rmk}

\begin{theo}{\rm (\cite{CR}  Theorem 3.10)}  Let $H$ be a finite
  dimensional Hopf algebra. Then there exists an associative algebra
  $X$ such that there is a vector space-preserving
  equivalence of categories between the category of left modules over
  $X$ and the category $\hopf$ of  Hopf bimodules over $H.$  \end{theo}

We can now state:

\begin{theo}{\rm (\cite{Ta}, \cite{Ta2})}\label{unificationtheorem} Let $H$ be a finite dimensional Hopf algebra. Then the
  following isomorphisms hold for any Hopf bimodules $M$ and $N$:

\begin{description} 

\item [(a)] $\Hh_{GS}^*(M,N) \cong \Ext_X^*(M,N)$

\item [(b)]$\Hh_{H4}^*(M,N) \cong \Ext_X^*(M,N)$ \end{description}
\end{theo}

\begin{proof} (This proof is similar to that in \cite{Ta} or
  \cite{Ta2}, but here we use injectives instead of projectives). Each
  isomorphism is proved using the universal property of $\Ext_X^*$; let us
remark that the category of Hopf bimodules over $H$ has enough
injectives, since it is equivalent to the category of modules over $X.$
  The functor $\Ext _{X}^*$ is then characterized by the following ({\it
cf.} \cite{McL}): \begin{description} \item [(1)] $ \Ext  _{X}^0 (M,N) \cong \Hom _{X}
(M,N) = \Hom _{\hopf} (M,N),$ 
\item [(2)] $ \Ext  _{X}^n (M,I) =0$ for every $n \geqslant 1$ and every
injective Hopf bimodule 
 $I,$ 
\item [(3)] $\Ext  _{X}^* (M,-) $ is a cohomological $\delta -$functor
(\emph{see}~\cite{W} p30).\end{description} Therefore we need to prove
that these three properties are satisfied for $\Hh_{GS}^*$ and $\Hh_{H4}^*$.

The proof of {\bf (1)} is straightforward in both cases, and the proofs of
{\bf (3)} are similar to the proofs given in \cite{Ta} and \cite{Ta2}, they
rely on alternative definitions of the double complexes and on
the structure of Hopf bimodules. We
shall now go through the proof of {\bf (2)} in each case.

\begin{lemma} \label{lemma2a} For every injective Hopf bimodule $I$ and
every integer $n \geqslant 1,$ the $k-$vector space $\Hh_{GS}^n(M,I)$
vanishes.\end{lemma}

\begin{prooflemma} Let $I$ be an injective Hopf bimodule. Then $0
  \rightarrow I
  \rightarrow I \rightarrow  0 \rightarrow  \cdots \rightarrow 0 \rightarrow \cdots $ is a relative  injective resolution of $I.$ It
is therefore homotopy equivalent to the cobar resolution $C^{\bullet}(I)$
(Proposition \ref{equivse}). Gerstenhaber and Schack's cohomology is
therefore the cohomology of the double complex in which all the terms
are zero except those on the first line, which are equal to
$\Hom_{\hopf}(B_q (M),I),$ for $q\geqslant 0.$ This line is acyclic, since $I$ is
injective and $B_\bullet (M)$ is exact. Its cohomology is therefore
zero in positive degree.  \end{prooflemma}

\begin{lemma}\label{lemma2b} For every injective Hopf bimodule  $I$ and every
integer $n \geqslant 1,$ the $k-$vector space $\Hh_{H4}^n(M,I)$ vanishes.
\end{lemma}

\begin{prooflemma} Let $I$ be an injective Hopf bimodule. Let us consider its cobar
resolution $$Cob^{\bullet}(I): 0 \rightarrow I
\stackrel{\lambda^{-1}}{\longrightarrow} H \oot I
\stackrel{\lambda^0}{\longrightarrow}   H^{\oot 2} \oot I
\rightarrow \cdots \rightarrow H^{\oot p+1} \oot I
\stackrel{\lambda^p}{\longrightarrow}  H^{\oot p+2} \oot I \rightarrow  \cdots,$$ with
$\lambda ^{-1}(u)=\rho _L (u)=u_{(-1)}\oot u_{(0)}.$

Since $I$ is injective and $\lambda ^{-1}$ is one-to-one, $\lambda ^{-1}$
has a retraction: there exists a morphism of Hopf bimodules $r: H \oot
I
\rightarrow  I$ satisfying $r(\lambda ^{-1}(u))=u$ for every $u
\in I.$

Set  \begin{eqnarray*} \chi^p\ :H^{\oot p+2}\oot I & \longrightarrow &
H^{\oot p+1}\oot I\\  h_0 \oot \ldots \oot h_{p+1} \oot u  &\mapsto &
h_0 \oot
 \ldots \oot h_p \oot r( h_{p+1} \oot u).\end{eqnarray*}
It is a Hopf bimodule morphism, and  $(\lambda ^{p-1}
\chi^{p-1} + \chi^{p}\lambda ^{p})= \id.$ Therefore $\chi^{\bullet} $ is a Hopf bimodule
homotopy from $\id$ to $0.$

Now fix $q \in \mathbb{N},$ and consider the complex $\Hom_{\hopf}(Bar_q(M),Cob^\bullet(I));$ the homotopy
$\chi^{\bullet}$ on $Cob^{\bullet}(I)$ yields a
homotopy $ \chi^{\bullet} \circ -$ from $\id$ to $0$ on this
complex. Therefore, the double complex: 
{\small $$\begin{array}{ccccccc}
\vdots & &\vdots & & \vdots &\\
\Hom_{\hopf}(Bar_1(M),I)
&\longrightarrow &\Hom_{\hopf}(Bar_1(M),Cob^0(I))
&\longrightarrow &\Hom_{\hopf}(Bar_1(M),Cob^1(I)) &\ldots  \\
\uparrow && \uparrow && \uparrow \\
\Hom_{\hopf}(Bar_0(M),I)
&\longrightarrow &\Hom_{\hopf}(Bar_0(M),Cob^0(I)) &\longrightarrow
&\Hom_{\hopf}(Bar_0(M),Cob^1(I)) &\ldots  \\ \\
\end{array}$$} which is the double complex defining $\Hh_{H4}^*(M,I)$ to which we have added one column, has exact rows. Its homology is
therefore that of the first column $\Hom_{\hopf}(Bar_\bullet(M),I)
$ (\emph{cf.} \cite{W} p59-60). Since $I$ is injective and $Bar_\bullet(M)$ is exact, the homology of this complex is zero in
positive degree.\end{prooflemma}

These two lemmas prove the second part of the characterization of
$\Ext_X^*$ for both cohomologies. \end{proof}

\begin{rmk}\label{xinvolvement} We have only used the existence of $X$
to say that the category of Hopf bimodules has enough injectives. \end{rmk}

\section{On the category of Hopf bimodules over an infinite
dimensional Hopf algebra; unification of cohomologies associated to an
infinite dimensional Hopf algebra}\label{enoughinjectives}

In this section, we shall prove that the category of Hopf bimodules
over any Hopf algebra has enough
injectives, and then extend the unification of Section
\ref{cohomologies}.

\begin{theo}\label{injectives} Let $H$ be any Hopf algebra. Then every
  Hopf bimodule over $H$  can be embedded in an injective
  Hopf bimodule. \end{theo}

\begin{proof} Let $M$ be a Hopf bimodule over $H$. Then $M$ can be
  viewed as an $H$-bimodule, and as such can be embedded in an
  injective $H$-bimodule $I$. Let $\varphi: M \hookrightarrow I$ denote
  this embedding. 

Now consider the Hopf bimodule $H\oot I \oot H$ (\emph{see} Remark 
\ref{inj}). We shall now prove that this is an injective Hopf bimodule
and that there is an embedding of Hopf bimodules from  $M$
to  $H\oot I \oot H$.

Let us consider the functor $\Hom_{\hopf}(-,H\oot I \oot H).$ Since
$\Hom_{\hopf}(-,H\oot I \oot H) $ is isomorphic to $ \Hom_{H-H}(-,I)$
(\emph{see} \cite{Ta2} Remark 1.22) and $I$ is an
injective bimodule, this functor is exact, thus proving that $H\oot I
\oot H$ is an injective Hopf bimodule.

Now define a map $\psi$ from $M$ to $H\oot I \oot H$ by setting
$\psi(m):=m_{(-1)}\oot \varphi(m_{(0)}) \oot m_{(1)}.$ This is a Hopf
bimodule map, and is clearly injective, since $(\varepsilon \oot \id
\oot \varepsilon) \circ \psi$ is equal to $\varphi$ which is injective.
\end{proof}

\begin{cor}\label{unifinfdim} Let $H$ be any Hopf algebra. Then the
  following isomorphisms hold for any Hopf bimodules $M$ and $N$:

\begin{description} 

\item [(a)] $\Hh_{GS}^*(M,N) \cong \Ext_\hopf^* (M,N)$

\item [(b)]$\Hh_{H4}^*(M,N) \cong \Ext_\hopf^* (M,N)$ \end{description}
\end{cor}

\begin{proof} This is proved in the same way as Theorem
  \ref{unificationtheorem}, now that we know that there are enough
  injective Hopf bimodules in the category $\hopf$. Note that we can consider the extensions $\Ext_{\hopf}^*(M,N)$, since the
category $\hopf$ is abelian.\end{proof}

\begin{rmk} Recall that when $H$ is finite dimensional, the category
  of Yetter-Drinfel'd modules is equivalent to the category of modules
over the Drinfel'd double $\mathcal{D}(H)$ (as braided categories), and therefore all these
cohomologies are isomorphic to $\Ext_{\mathcal{D}(H)}^*(M^R,N^R).$ 

In
particular, $\hb \cong \Ext_{\mathcal{D}(H)}^* (k,k)$, which answers a
question of F. Panaite and D. \c Stefan (\emph{see} \cite{PS}). \end{rmk}

\begin{rmk}  We proved in \cite{Ta} (Corollary 2.18) that the
  cohomologies we have considered here are Morita invariant when the
  Hopf algebra is finite dimensional. It is clear that this result
  extends to the infinite dimensional case, that is
  that if $H$ and $H'$ are Hopf algebras which are Morita equivalent
  as algebras and such that the functors $\mathcal{F}$ and
  $\mathcal{G}$ giving the equivalence are
  monoidal, then there is an isomorphism $$\Ext_{\hopf}^*(M,N) \cong
  \Ext_{\hopf'}^*(\mathcal{F}(M),\mathcal{F}(N))$$ for all $H$-Hopf
  bimodules  $M$ and $N$. Here, $\hopf'$ denotes the category of Hopf
  bimodules over $H'$. In particular,  $$\Hh _b^*(H,H) \cong \Hh _b^*(H',H').$$ 
\end{rmk}

\section{Graded-commutativity of the Yoneda product of Hopf \\bimodule extensions}\label{commut}

In \cite{Ta} and \cite{Ta2}, we have defined a cup-product on $\Hh_{H4}^*$
which corresponds to the Yoneda product of extensions up to sign \emph{via} the
identification of Corollary \ref{unifinfdim}
 (the proof of the
correspondence of the products remains valid for an infinite
dimensional Hopf algebra). This part is devoted to the  study of this product,
in particular proving that it is graded-commutative, using techniques
of S. Schwede \cite{S}.

\subsection{The products on the cohomologies $\Hh_{H4}^*$ and $\Hh_b^*$}

Let us first define the cup-product on the cohomology $\Hh_{H4}^*:$

\begin{prop}{\rm (\cite{Ta}, \cite{Ta2})} Let $f \in \Hom_{H-}^{-H} (M \otimes H ^{\otimes p-s} ,H
^{\otimes s} \otimes L)$ be a $p-$cochain (a map of left $H$-modules
and right $H$-comodules) and $g \in \Hom_{H-}^{-H} (L
\otimes H ^{\otimes q-r} ,H ^{\otimes r} \otimes N)$ be a
$q-$cochain. Set $n=p+q$ and $t=s+r.$ Define the $n-$cochain $f \smile
g \in \Hom_{H-}^{-H} (M \otimes H ^{\otimes n-t} ,H ^{\otimes t}
\otimes N)$ by:  \begin{align*} f \smallsmile g &(m \otimes \mathbf{a_{1, n-t}} )=\\
& (-1)^{s(q-r)}(1 ^{\otimes s} \otimes g )\,[f(m \otimes
\mathbf{a_{1,p-s}} ) .(\Delta ^{(s-1)} (a^{(1)}_{p-s+1} \ldots a^{(1)}_{n-t} ) \otimes 1 ) 
 \otimes  \mathbf{a^{(2)}_{p-s+1,n-t}}].
\end{align*}

The differential $D$ of the total complex associated to the
Hopf bimodule double complex is a right derivation for the cup-product
$\smile,$ that is $$D(f \smile g) = Df \smile g + (-1)^p f \smile
Dg,$$ so that the formula for $\smallsmile$ yields a product $\Hh_{H4}^* (M,L) \otimes \Hh_{H4}^* (L,N) \rightarrow \Hh_{H4}^* (M,N). $
\end{prop}

\begin{rmk}   The cup-product in $\Hh _b^*(H,H)$ is as follows:
let $f$ be in  $\Hom _k (H^{\otimes p-s} , H^{\otimes s})$ and $g$
in $ \Hom _k (H^{\otimes q-r}, H^{\otimes r});$ then   \begin{eqnarray*}
(f \smallsmile g) (a_1 \otimes \ldots \otimes a_{n-t}) & = & (-1)^{s(q-r)} f(a^{(1)}_1 \otimes \ldots \otimes a^{(1)}_{p-s} ) \ \Delta ^{(s-1)} (a^{(1)}_{p-s+1} \ldots a^{(1)}_{n-t}) \\
& &\  \otimes \Delta ^{(r-1)} (a^{(2)}_{1} \ldots a^{(2)}_{p-s})\
g(a^{(2)}_{p-s+1} \otimes \ldots \otimes a^{(2)}_{n-t}).
\end{eqnarray*} 
 \end{rmk}

We shall now relate the cup-product $\smallsmile$ with the Yoneda
product of extensions.

\begin{theo}{\rm (\cite{Ta}, \cite{Ta2})}\label{Yonedathe} Let $H$ be a Hopf
algebra. Let $M, \ N$ and $L$ be Hopf bimodules
over $H$. Let $\varphi _{MN} : \Hh_{H4} ^* (M,N) \rightarrow
\Ext_{X}^* (M,N)$ be the isomorphism  extending $\id_{\Hom_{\hopf}
  (M,N)}$ and let $\sharp$ denote the Yoneda product of extensions
(\emph{see} below).

Then, if  $f \in \Hh _{H4}^p(M,L)$ and $g \in \Hh _{H4}^q (L,N),$ the relationship between the products is given by $$\varphi _{MN}^{p+q} (f \smile g) = (-1)^{pq} \, \varphi _{LN}^q (g) \sharp \varphi _{ML}^p (f):$$ the cup product and the Yoneda product are equal up to sign.  
\end{theo}

\begin{proof}(Sketch) After proving universal properties of the cup-product and
the Yoneda product, we reduce to the case when $q=0;$ this case is
proved by induction on $p,$ thanks to a partial associativity
property, and the knowledge of $\varphi ^1$.
 \end{proof}

Therefore, to study the algebraic structure of $\Hh_b^*(H,H)$, we may
consider the algebra $\Ext_{\hopf}^*(H,H)$ endowed with the Yoneda
product.


\subsection{Graded-commutativity of the Yoneda product}

We will follow S. Schwede's method (\cite{S}) to prove that the Yoneda
product of Hopf bimodule extensions of $H$ by $H$
 is graded-commutative. We refer to \cite{S} for the definitions and
 results on the homotopy of (the nerve of) a category which we shall need. Let us fix notations.

We shall denote by $\ext ^n (H,H)$ the category of extensions of Hopf
bimodules from $H$ to $H$. Then $\pi _0 \ext ^n (H,H) \cong
\Ext_{\hopf}^n(H,H)$ and $\pi _1 \ext ^n (H,H) \cong \Ext_{\hopf}
^{n-1}(H,H),$ where paths and loops are defined as in \cite{S}. Now
given two extensions  ${\bf E}:0\stackrel{i_{\e}}{\rightarrow} H \rightarrow E_{m-1}\rightarrow
  E_{m-2}\rightarrow \cdots \rightarrow
  E_0\stackrel{p_{\e}}{\rightarrow} H\rightarrow 0$ in $\ext^m(H,H)$ and ${\bf F}:0\rightarrow H \rightarrow F_{n-1}\rightarrow
  F_{n-2}\rightarrow \cdots \rightarrow F_0\rightarrow
  H\stackrel{p_{\f}}{\rightarrow} 0$  in $\ext^n(H,H),$ we define
  their Yoneda product to be the class in $\pi _0 \ext^{m+n}(H,H)$ of the
  extension obtained by `splicing' $\e$ and $\f$:
  $$\e \sharp \f: 0\rightarrow H\rightarrow F_{n-1}\rightarrow\cdots\rightarrow F_0
  \stackrel{i_{\e}p_{\f}}{\longrightarrow}E_{m-1}\rightarrow \cdots
  \rightarrow E_0 \rightarrow H\rightarrow 0.$$ We denote by $(-1)\e$
  the sequence obtained from $\e$ by replacing $p_{\e}$ by $-p_{\e};$
  this represents the inverse of $\e$ in $\pi _0 \ext^m (H,H)$ with
  respect to Baer sum. 

To prove that this product is graded-commutative, we need to prove
that the classes in $\pi _0 \ext^{m+n} (H,H)$ of $\f \sharp \e$ and of $(-1)^{mn}\e \sharp \f$ are
equal, so we need to find a path from  $\f \sharp \e$ to $(-1)^{mn}\e
\sharp \f$ in $\ext^{m+n}(H,H).$ This will involve a tensor product
of extensions of Hopf bimodules;
 the first step is
to define  this tensor product over $H$ of Hopf bimodules and
hence of extensions of Hopf bimodules, and to check that it makes sense in the category of
Hopf bimodules $\hopf.$

\begin{dfn} Let $E$ and $F$ be Hopf bimodules over $H.$ Define $E \qot
   F$ as the quotient $$\frac{E \uot F}{<\{eh\uot f -e \uot hf \mid \
    e \in E, \ h \in H, \ f \in F\}>}.$$ \end{dfn}

\begin{lemma} The space $E \qot F$ is a Hopf bimodule. \end{lemma}

\begin{proof} Let $I$ be the vector space generated by $\{eh\uot f -e \uot hf \mid \
    e \in E, \ h \in H, \ f \in F\}.$ It is well-known that $I$ is a
    sub-bimodule of $E \uot F.$ It is straightforward to check that it is also a
    sub-bicomodule of $E \uot F$ (recall that the coactions on $E
    \uot F$ are codiagonal). Therefore $I$ is a sub-Hopf bimodule of
    $E \uot F$, and the quotient is then a Hopf bimodule. \end{proof}

Now we need to see that the tensor product of extensions of Hopf
bimodules is again an extension; let us consider a special case first:

\begin{lemma}\label{tensorextension} If $E$ is a Hopf bimodule and $F' \rightarrow F
  \rightarrow F''$ is any exact sequence of Hopf bimodules, then
  $$(S)\ \ E \qot F'\rightarrow E \qot F\rightarrow E \qot
  F''$$ is an exact sequence of Hopf bimodules.\end{lemma}

\begin{proof} The sequence $(S)$ is a sequence of Hopf
  bimodules, since the spaces are Hopf bimodules by the previous lemma
  and the maps are tensor products of Hopf bimodule morphisms and are
  therefore Hopf bimodule morphisms.

It remains to be checked that the sequence is exact. Since $E$ is a
Hopf bimodule, we know that it is free as a
  right $H$-module (\emph{see} for instance
\cite{Mo} Theorem 1.9.4), therefore $E$ is flat as a right
$H$-module. So $(S)$ is exact as a sequence of vector spaces, hence
also as a sequence of Hopf bimodules.  \end{proof}
 
We may now consider tensor products of extensions:

\begin{prop} Let ${\bf E}:0\rightarrow H \rightarrow E_{m-1}\rightarrow
  E_{m-2}\rightarrow \cdots \rightarrow E_0\stackrel{p_{\e}}{\rightarrow} H\rightarrow 0$ and ${\bf F}:0\rightarrow H \rightarrow F_{n-1}\rightarrow
  F_{n-2}\rightarrow \cdots \rightarrow F_0\rightarrow H\stackrel{p_{\f}}{\rightarrow} 0$ be extensions of Hopf
  bimodules. Then their tensor product over $H, $ defined by $$({\bf E}
  \qot {\bf F})_r:=\bigoplus_{\scriptsize \begin{array}{c}s+t=r\\s,t\geqslant
      0\end{array}}E_s \qot F_t$$ for $0\leqslant r \leqslant m+n$ and $({\bf E}
  \qot {\bf F})_{-1}:=H$, is also a Hopf bimodule extension of $H$
  by itself.\end{prop}

\begin{proof} Fix $-1\leqslant s\leqslant m$: from Lemma \ref{tensorextension},
  the sequence $E_s \qot {\bf F}$ is an extension of Hopf
  bimodules. We then take the direct sum over $-1\leqslant s\leqslant m$ of
  these exact sequences of Hopf bimodules to obtain an extension of
  Hopf bimodules ${\bf E} \qot {\bf F}.$ \end{proof}

Given two extensions ${\bf E}$ and ${\bf F}$ as above, we may now construct a
path from ${\bf F} \sharp {\bf E}$ to $(-1)^{mn} {\bf E} \sharp {\bf
  F}.$ In fact, as in \cite{S},
 we are going to construct a loop in $\ext^{m+n}(H,H) $ going
 through these two extensions.

Consider the following two maps, $\lambda _{\e,\f}: \e \qot \f
\rightarrow \f \sharp \e$ and $\rho _{\e,\f}: \e \qot \f \rightarrow
(-1)^{mn}\e \sharp \f,$ defined on each degree by:

\begin{eqnarray*}
\mbox{if } 0 \leqslant i < m,&( \lambda_{\e,\f}) _i:& (\e \qot \f)_i
\stackrel{projection}{\longrightarrow} E_0 \qot F_i
\stackrel{p_{\e}\qot \id _{F_i}}{\longrightarrow} F_i \\
\mbox{if } m \leqslant i \leqslant m+n,& ( \lambda_{\e,\f}) _i:&(\e \qot \f)_i
\stackrel{projection}{\longrightarrow} E_{i-m} \qot F_m
\stackrel{\sim}{\longrightarrow} E_{i-m} \\
\mbox{if } 0 \leqslant j < n,& (\rho_{\e,\f})_j :& (-1)^{mn} \left[
(\e \qot \f)_j \stackrel{projection}{\longrightarrow} E_j \qot F_0
\stackrel{\id _{\e _ j} \qot p_{\f}}{\longrightarrow} E_j \right]\\
\mbox{if } n \leqslant j \leqslant m+n,& (\rho_{\e , \f})_j :& (-1)^{m+n-j} \left[(\e \qot \f)_j
\stackrel{projection}{\longrightarrow} E_n \qot F_{j-n}
\stackrel{\sim}{\longrightarrow} F_{j-n}\right]. \end{eqnarray*}

These are morphims of complexes (this is a straightforward
computation) of Hopf bimodules, since each of the maps is given by
composition of projections, natural identifications $E_i \qot H \cong
E_i$ or $H \qot F_j \cong F_j,$ and the tensor product of an identity map with a
morphism of Hopf bimodules, so they are all morphisms of Hopf
bimodules.

We may now consider the following loop of extensions in $\ext^{m+n} (H,H)$:

$$\xymatrix{& \e \qot \f \ar_{\lambda _{\e,\f}}[dl] \ar^{\rho
    _{\e,\f}}[dr]\\ 
\f \sharp \e && (-1)^{mn} \e \sharp \f \\  & (-1)^{mn} \f \qot \e
\ar ^{\rho _{\f,\e}}[ul] \ar _{\lambda _{\f,\e}} [ur]} $$ oriented counter-clockwise. We call it $\Omega(\f,\e).$

There is therefore a path in $\ext^{m+n} (H,H)$ from $\f
\sharp \e$ to $(-1)^{mn}\e \sharp \f$ (either the upper part or the
lower part of the diagram above), so  $\f
\sharp \e$ and $(-1)^{mn}\e \sharp \f$ define the same element in
$\pi_0 \ext^{m+n} (H,H)\cong \Ext_{\hopf}^{m+n}(H,H).$

Furthermore, as in \cite{S}, the construction of the loop $\Omega(\f,\e)$ is functorial in $\f$ and $\e,$ so
this defines a map: $$\Omega: \pi_0 \ext^m (H,H)\times
\pi_0\ext^n (H,H) \longrightarrow \pi_1
\ext^{m+n} (H,H),$$ called the loop bracket.

\begin{rmk} The graded-commutativity can be deduced from the above
  construction without involving homotopy groups: indeed, we have
  constructed maps of extensions between $\e \qot \f $ and $\f \sharp
  \e$, and between $\e \qot \f$ and  $ (-1)^{mn} \e \sharp
  \f$. Therefore, by definition of the equivalence relation between
  extensions (\emph{see} \cite{McL}), these three extensions are
  equivalent, and represent the same element in
  $\Ext^{m+n}_{\hopf}(H,H)$. 

The homotopy groups become useful when defining the bracket, which we
shall now discuss further. \end{rmk}

\section{The loop bracket}\label{bracket}

The cohomologies we have described above ($\Hh_{A4}^*$ and $\Hh_b^*$) have many points in common
with Hochschild cohomology of algebras. They involve bar resolutions,
are isomorphic to an $\Ext^*$ functor, and are endowed with a 
cup-product which corresponds to the Yoneda product of extensions  and which is
graded-commutative. It is therefore natural to wonder if $\Hh_b^*(H,H)$
might be a G-algebra, that is if it might be endowed with a graded Lie
bracket compatible with the cup-product, as M. Gerstenhaber proved is
the  case for Hochschild cohomology. We have been unable to answer
this question, but we describe below some points in favour of this
being the case. 

Let us recall the definition of a G-algebra:

\begin{dfn} A \emph{G-algebra} is a graded $k$-module $\Lambda
  = \oplus _n \Lambda_n$ equipped with two multiplications, $(\lambda,
  \nu) \mapsto \lambda \smallsmile \nu$ and $ (\lambda,\nu)\mapsto
  [\lambda,\nu],$ satisfying the following properties:
\begin{description}
\item [(1)] $\smallsmile$ is an associative graded (by degree)
  commutative product;
\item [(2)] $[-,-]$ is a graded Lie bracket for which the grading is
  reduced degree, this being one less than the degree;
\item [(3)] $[-,\eta ^p]$ is a degree $p-1$ graded derivation of the
  associative algebra structure for all $\eta ^p \in \Lambda _p:$
  $$[\lambda ^m \smallsmile \nu ^n,\eta ^p]=[\lambda,\eta]\smallsmile
  \nu + (-1)^{m(p-1)} \lambda \smallsmile [\nu,\eta] \quad \forall
  \lambda, \nu.$$
\end{description}
\end{dfn}

We know that the first condition is satisfied (Section
\ref{commut}). Regarding the bracket, in the previous section we described a map $$\Omega: \pi_0 \ext^m (H,H)\times
\pi_0\ext^n (H,H) \longrightarrow \pi_1
\ext^{m+n} (H,H).$$ Transporting this \emph{via} the isomorphisms $\varphi$
between $\Hh_b^{*}(H,H)$ and $\pi _0 \ext ^{*}(H,H)$ and $\psi$ between
$\Hh_b^{*-1}(H,H)$ and $\pi _1 \ext ^{*}(H,H)$ gives maps
$\Hh_b^m(H,H) \times \Hh_b^n(H,H)\rightarrow \Hh_b^{m+n-1}(H,H)$. We
therefore have a candidate for the product $[-,-].$
 However, we do not know any of the isomorphisms explicitly, so we
 cannot prove directly conditions {\bf (2)} and {\bf (3)}. More precisely, since
 we want to keep as close as possible to the Hochschild situation
 (\emph{see} \cite{S}), we consider the map $[-,-]: \Hh_b^m(H,H) \times
 \Hh_b^n(H,H)\rightarrow \Hh_b^{m+n-1}(H,H)$ such that the diagram
 $$\xymatrix{\Hh_b^m(H,H) \times \Hh_b^n(H,H) \ar^(.57){[-,-]}[r]\ar_{\varphi\times\varphi}[d]&
   \Hh_b^{m+n-1}(H,H)\ar^{\psi}[d]\\\pi_0\ext^m(H,H)\times
   \pi_0\ext^n(H,H)\ar_(.62){\Omega}[r]&\pi_1\ext^{m+n}(H,H) }$$ commutes
 \emph{up to the sign $(-1)^n$}, that is, if  $\e=\varphi(e)$ is of
 degree $m$ and $\f=\varphi(f)$ is of degree $n$, then
 $\Omega(\e,\f)=\psi((-1)^n[e,f])$.
 However, we do not know any of the
 isomorphisms explicitly, so we do not have an expression for $[-,-]$,
 and we
 cannot prove directly conditions {\bf (2)} and {\bf (3)}.

Therefore, we now need to find out how the conditions {\bf (2)} and {\bf (3)} translate for
$\Omega$. These conditions involve in particular the sum of elements of
$\Hh_b^*(H,H)$. This  corresponds to concatenation of loops:

\begin{prop} The isomorphism $\psi$ is an isomorphism of
  groups, where the group law on $\Hh_b^*(H,H)$ is induced by Baer sum
and that on $\pi_1\ext^{*+1}(H,H)$ is induced by concatenation of
loops. \end{prop}

\begin{proof} We have an isomorphism $\Hh_b^*(H,H)
  \stackrel{\sim}{\rightarrow} \pi _0 \ext^*(H,H)=\Ext_{\hopf}^*(H,H)$: this is the
  isomorphism of Corollary \ref{unifinfdim},
 which is a morphism of
  (abelian) groups (\emph{cf.} \cite{McL}), where the group law on extensions is Baer sum. We are interested in the
  isomorphism    $\Hh_b^*(H,H)
  \stackrel{\sim}{\rightarrow} \pi _1
  \ext^{*+1}(H,H)=\Ext_{\hopf}^*(H,H)$ which we know exists from
  \cite{R}. We know (\cite{R} Theorem 2) that the isomorphism $\pi _0
  \ext ^n (H,H) \cong \pi _1
  \ext ^{n+1} (H,H) $ comes from a  functor $F_{n+1}$ which is a
  \emph{group} isomorphism (with Baer sum on the left and
  concatenation of loops on the right). Therefore by composition, the
  isomorphism $\psi$ is an isomorphism of groups. \end{proof}

We now wish to prove  properties {\bf (2)} and {\bf (3)} adapted to
$\Omega.$ They are satisfied in the case of an associative algebra
when considering bimodule extensions, since S. Schwede proved in \cite{S} that
$\Omega$ corresponds to Gerstenhaber's Lie bracket in that case.
However, he doesn't give a direct proof of these properties for
$\Omega,$ and describing the Jacobi identity for $\Omega$ requires the
construction of a loop (or of an element in a higher homotopy group) associated not to two extensions but to one
extension and one loop, and property {\bf (3)} requires the definition of the
Yoneda product of a loop and an extension. We do not yet know how to
do this;  however, the graded
anti-commutativity of the bracket can easily be translated as follows:

\begin{prop} The graded anti-commutativity of the Lie bracket
  translates as $$\Omega(\e,\f)=\Omega(\f,\e)^{(-1)^{mn}}$$ where $\e$ is of
degree $m$ and $\f$ is of degree $n$. \end{prop}

\begin{proof} Indeed, if $\e=\varphi(e)$ and $\f=\varphi(f)$, we have
  $$\begin{array}{rcl}\Omega(\e,\f)&=&\psi((-1)^n[e,f])\\
&=&\psi((-1)^n(-1)^{(m-1)(n-1)+1}[f,e])\\
&=&\psi((-1)^m[f,e])^{(-1)^{mn}}\\
&=&\Omega(\f,\e)^{(-1)^{mn}}.\end{array}$$ \end{proof}

Let us now prove that $\Omega$ does indeed satisfy this property:

\begin{prop}\label{antisymmetry} The loop bracket $\Omega$ satisfies:  $$\Omega(\e,\f)=\Omega(\f,\e)^{(-1)^{mn}}$$ where $\e$ is of
degree $m$ and $\f$ is of degree $n$. \end{prop}

\begin{proof} The first case we consider is when either $m$  or $n$ is
even. Then we wish to prove that $\Omega(\e,\f)=\Omega(\f,\e)$, and
this is obvious when looking at the loops: one is obtained from the
other by rotating by 180 degrees. 

In fact, we notice that in general
$\Omega(\e,\f)=\Omega(\f,(-1)^{mn}\e)$. Therefore in the case when
both $m$ and $n$ are odd, we wish to prove that
$\Omega(\f,-\e)=\Omega(\f,\e)^{-1},$ or that
$\Omega(\f,\e).\Omega(\f,-\e)$ is the trivial loop. To do this, we
shall consider  the behaviour of $\Omega$ towards Baer sums.

Note that the components of the category of extensions have
canonically isomorphic groups, the isomorphisms are given by Baer sum
with a fixed extension (and converses are given by Baer sum with the
opposite extension). Therefore, even if the loops we are going to
consider are not always in the same component originally, they may be
considered so.

Now given three extensions $\e,\ \f,$ and $\g,$ we wish to compare the
loops
$\Omega (\g, \e+\f)$ and $\Omega(\g,\e).\Omega(\g,\f)$. Since both the
Yoneda product $\sharp$ and the tensor product $\qot$ commute with
Baer sum, the left hand term is obtained by adding the extensions in
the loops $\Omega(\g,\e)$ and $\Omega(\g,\f)$ componentwise. 

Therefore the question may be set more generally as follows: given two
loops of extensions $L: Y \stackrel{\alpha _0}{\longleftrightarrow}
C_1  \stackrel{\alpha _1}{\longleftrightarrow} \cdots \stackrel{\alpha
  _{r-1}}{\longleftrightarrow} C_r  \stackrel{\alpha
  _r}{\longleftrightarrow} Y $ and $L': Y' \stackrel{\alpha' _0}{\longleftrightarrow}
C'_1  \stackrel{\alpha' _1}{\longleftrightarrow} \cdots \stackrel{\alpha'
  _{r-1}}{\longleftrightarrow} C'_r  \stackrel{\alpha'
  _r}{\longleftrightarrow} Y' $, we wish to compare the loops of
extensions $L'': Y+Y' \stackrel{(\alpha _0,\alpha' _0)}{\longleftrightarrow}
C_1+C'_1  \stackrel{(\alpha _1,\alpha' _1)}{\longleftrightarrow} \cdots \stackrel{(\alpha
  _{r-1},\alpha'
  _{r-1})}{\longleftrightarrow} C_r + C'_r \stackrel{(\alpha
  _r,\alpha'
  _r)}{\longleftrightarrow} Y+Y' $ and $L''': Y +Y'\stackrel{(\alpha _0,\id)}{\longleftrightarrow}
C_1 +Y' \stackrel{(\alpha _1,\id)}{\longleftrightarrow} \cdots \stackrel{(\alpha
  _{r-1},\id)}{\longleftrightarrow} C_r+Y'  \stackrel{(\alpha
  _r,\id)}{\longleftrightarrow} Y+Y' \stackrel{(\id,\alpha' _0)}{\longleftrightarrow}
Y +C'_1  \stackrel{(\id,\alpha' _1)}{\longleftrightarrow} \cdots \stackrel{(\id,\alpha'
  _{r-1})}{\longleftrightarrow} Y+ C'_r  \stackrel{(\id,\alpha'
  _r)}{\longleftrightarrow} Y+Y' $. Note that in the last loop, the
original loops $L$ and $L'$ have been shifted so that they are in the
same component, based at $Y+Y',$ and can therefore be composed.

\begin{lemma} The loops of extensions $L''$ and $L'''$ are homotopic.
\end{lemma}

\begin{prooflemma} To simplify notation, we shall look at the case when
  r=1; the general case is similar.

In this case, $L''$ takes the form $Y+Y' \stackrel{(\alpha _0,\alpha' _0 )}{\longleftrightarrow}
C_1+C'_1  \stackrel{(\alpha _1,\alpha' _1)}{\longleftrightarrow}
Y+Y',$ which is homotopic to $Y+Y'\stackrel{(\alpha
  _0,\id)}{\longleftrightarrow}C_1 + Y' \stackrel{(\id,\alpha'
  _0)}{\longleftrightarrow} C_1+C'_1 \stackrel{(\alpha
  _1,\id)}{\longleftrightarrow}Y+C'_1\stackrel{(\id,\alpha'
  _1)}{\longleftrightarrow}Y+Y'$ by definition of the homotopy of paths. This loop is homotopic to $Y+Y'\stackrel{(\alpha
  _0,\id)}{\longleftrightarrow}C_1 + Y' \stackrel{(\alpha_1,\alpha'
  _0)}{\longleftrightarrow}Y+C'_1 \stackrel{(\id,\alpha'
  _1)}{\longleftrightarrow}Y+Y'$ which is finally homotopic to $Y+Y'\stackrel{(\alpha
  _0,\id)}{\longleftrightarrow}C_1 + Y'    \stackrel{(\alpha
  _1,\id)}{\longleftrightarrow} Y+Y'   \stackrel{(\id,\alpha'
  _0)}{\longleftrightarrow}              Y+C'_1 \stackrel{(\id,\alpha'
  _1)}{\longleftrightarrow}Y+Y'=L'''.  $\end{prooflemma}

Finally we obtain:

\begin{lemma}\label{bracketbaersum} Given three extensions $\e, \ \f,$ and $\g,$ and changing
base points so that the identities make sense, we have:
\begin{eqnarray*} \Omega(\g,\e+\f)&=&\Omega(\g,\e). \Omega(\g,\f)
  \mbox{ and}\\ 
\Omega(\f+\g,\e)&=&\Omega(\f,\e). \Omega(\g,\e).
\end{eqnarray*}\end{lemma}

Now using Lemma \ref{bracketbaersum}, we know that
$\Omega(\f,\e).\Omega(\f,-\e)=\Omega(\f,{\bf 0}),$ where ${\bf 0}$ is
the trivial extension of degree $m.$ So if
$\Omega(\f,{\bf 0})$ is the trivial loop, we have proved Proposition \ref{antisymmetry}.

Let us consider the loop $\Omega(\f,{\bf 0})$: computing the maps in
this loop gives $\lambda_{{\bf 0},\f}=\rho _{\f,{\bf 0}}$ and $\lambda
_{\f,{\bf 0}}= \rho _{{\bf 0},\f}$. We also have $\f \qot {\bf 0}=
{\bf 0} \qot \f$ and $\f \sharp  {\bf 0}= {\bf 0} \sharp
\f$. Therefore, using the homotopy relations, the loops 
$$\xymatrix{&{\bf 0} \qot \f  \ar_{\lambda _{{\bf 0},\f}}[dl] \ar^{
   \lambda _{\f,{\bf 0}}}[dr]&&&{\bf 0} \qot \f   \ar^{
   \lambda _{\f,{\bf 0}}}[d]\\ \f \sharp {\bf 0}& \Omega(\f,{\bf 0})& {\bf 0} \sharp \f&\mbox{and}& {\bf
  0}\sharp \f \\& \f \qot {\bf 0}\ar^{\lambda _{{\bf 0},\f}}[ul]
  \ar_{\lambda _{\f,{\bf 0}}}[ur]&&&{\bf 0} \qot \f  \ar_{\lambda _{\f,{\bf 0}}}[u]}$$
are homotopic, and
the last one is homotopic to the trivial loop based at ${\bf
  0} \qot \f $. This concludes the proof.
\end{proof}

\end{document}